\documentclass[11pt,reqno]{amsart}
\usepackage{amsmath,amssymb,amsxtra,latexsym,amsthm,amscd,amsfonts,mathrsfs,pb-diagram}
\usepackage[colorlinks=true,linkcolor=magenta,citecolor=blue]{hyperref}
\usepackage[margin=1in]{geometry}
\usepackage{enumitem}
\usepackage{mathabx}
\usepackage{epsfig}

\usepackage{color}
\def\Int{\displaystyle\int}

\begin{document}
	\title[\hfil Generalized damped Klein-Gordon equation with nonlinear memory]{on generalized damped Klein-Gordon equation with nonlinear memory}
\thanks{$\star$ Corresponding author}
	{}  
	{}  
\subjclass[2010]{35A01, 35B45, 35G30}

\keywords{Klein-Gordon equation; frictional damping; mass term; global existence; nonlinear memory.}

	\maketitle
	\centerline{\scshape Said Khaldi$^*$}
	\medskip
	
	{\footnotesize
	\centerline{Laboratory of Analysis and Control of PDEs, Djillali Liabes University}
	\centerline{P.O. Box 89, Sidi-Bel-Abbes 22000, Algeria}
	\centerline{Emails: saidookhaldi@gmail.com, said.khaldi@univ-sba.dz}
	}
\medskip

\centerline{\scshape Mohamed Menad}
\medskip
{\footnotesize
	\centerline{Laboratoire de Mathematiques et Applications (LMA), Universitié Hassiba Benbouali de Chlef} 
	\centerline{Hay Essalem B.P. 151,
		02000, Chlef, Algerie}
	\centerline{Email: menmo2001@gmail.com}
} 

\begin{abstract}
In this paper we consider the Cauchy problem for linear dissipative generalized Klein-Gordon equations with nonlinear memory in the right hand side. Our goal is to study the effect of this nonlinearity on both the decay estimates of global solutions as well as the admissible range of the exponent $p$.
\end{abstract}
\numberwithin{equation}{section}
\newtheorem{theorem}{Theorem}[section]
\newtheorem{lemma}[theorem]{Lemma}
\newtheorem{definition}[theorem]{Definition}
\newtheorem{remark}[theorem]{Remark}
\newtheorem{proposition}[theorem]{Proposition}
\allowdisplaybreaks		
\section{Introduction}\label{sec:section1}
We are concerned with the Cauchy problem for the linear dissipative generalized Klein-Gordon equations with nonlinear memory in the right hand side:
\begin{equation}\label{1.1}
\left\lbrace 
\begin{array}{ll}
\partial_{t}^{2}u+(-\Delta)^{\sigma}u+2a\partial_{t}u+m^{2}u=\Int_{0}^{t}(t-s)^{-\gamma}\left|u(s,x)\right| ^{p}ds\\
\\
u(0,x)=u_{0}(x),   \ \ \ \partial_{t}u(0,x)=u_{1}(x)
\end{array}\right. (t,x)\in \mathbb{R_{+}}\times \mathbb{R}^{n},
\end{equation}
where 
\begin{equation*}
a>0, \ \ m>0, \ \ \gamma\in(0,1), \ \ p>1 \ \ \ and \ \ \sigma\geq 1.
\end{equation*}
Usually, the terms $2a\partial_{t}u$ and $m^{2}u$ are called frictional damping and mass, respectively. The operator $(-\Delta)^{\sigma}$ denotes the fractional Laplacian on $\mathbb{R}^{n}$ with symbol $|\xi|^{2\sigma}$, i.e, is defined through the Fourier transform $\mathcal{F}$ as follows : 
$$
\mathcal{F}\left( (-\Delta)^{\sigma}f\right) = |\xi|^{2\sigma}\mathcal{F}\left(f\right) (\xi), \ \  \xi\in\mathbb{R}^{n}, \ \  |\xi|^{2}=\left( \sum_{i=1}^{n}\xi_{i}^{2}\right).
$$  
Let us recall a previous result about the Cauchy problem (\ref{1.1}) but with a power nonlinearity in the right hand side, that is,
\begin{equation}\label{1.2}
\left\lbrace 
\begin{array}{ll}
\partial_{t}^{2}u+(-\Delta)^{\sigma}u+2a\partial_{t}u+m^{2}u=|u|^{p}\\
\\
u(0,x)=u_{0}(x),   \ \ \ \partial_{t}u(0,x)=u_{1}(x)
\end{array}\right. (t,x)\in \mathbb{R_{+}}\times \mathbb{R}^{n}.
\end{equation}
This later problem has been successfully studied in \cite{d'abbicco2015} where the author have first derived some exponential decay in the $L^{2}-L^{2}$ estimates for solutions to the corresponding linear equations, namely,
\begin{equation}\label{1.3}
\left\lbrace 
\begin{array}{ll}
\partial_{t}^{2}u+(-\Delta)^{\sigma}u+2a\partial_{t}u+m^{2}u=0\\
\\
u(0,x)=u_{0}(x),   \ \ \ \partial_{t}u(0,x)=u_{1}(x)
\end{array}\right. (t,x)\in \mathbb{R_{+}}\times \mathbb{R}^{n},
\end{equation}  
and then proved the global (in time) existence of small data solution for any exponent $p$ satisfies: 
$$1<p \leq \frac{n}{n-2\sigma} \ \ \ if \ \ n>2\sigma. $$
He obtained the same exponential decay estimates as those of (\ref{1.3}). So the fact that power nonlinearity does not influence the decay estimates of solution to the semi-linear Cauchy problem is well known from this study \cite{d'abbicco2015} and other in references \cite{D'AbbiccoEbert}, \cite{EbertReissig}, \cite{PhamKainaneReissig}.\\
Now for the nonlinear memory (see for instance \cite{D'Abbiccomemory}) the situation is completely different, there appears some loss of decay for both the solution and related quantities like energy. Moreover, the so-called critical exponent also changes. This influence of a nonlinear memory motivates us to study the interesting Cauchy problem (\ref{1.3}) by putting it in the right hand side. \\
In this paper we also show a remarkable effect of this nonlinearity on the decay estimates of global solution, in fact, it will be satisfies some polynomial decay instead of exponential decay. Moreover, the exponent $p$ in (\ref{1.1}) will be scaled by the factor $\frac{1}{\gamma}$ with $\gamma \in (0,1)$.
\\
The present paper is organized as follows : in Section \ref{preliminaries and tools} we recall some important tools such as $L^{2}-L^{2}$ estimates and integral inequalities. Section \ref{GlobalExistence} contains the main result of global (in time) existence and its proof which is essentially based the Banach fixed point theorem.
\section{preliminaries and tools }\label{preliminaries and tools}
In this section we gather the tools that will be used to prove our result. First we recall the quite important $L^{2}-L^{2}$ estimates which is proved in \cite{d'abbicco2015}.
\begin{lemma}\cite[Proposition 2.1]{d'abbicco2015}\label{first lemma}
	Let $a>0$, $m>0$ and $\sigma>0$ in (\ref{1.3}). Then, for any real $k \geq 0$ and for any integer $j\geq  0$, the solution to (\ref{1.3}) satisfies the following decay estimates:
	\begin{align} 
	\|\partial_{t}^{j}u^{lin}(t,\cdot)\|_{H^{2k}} &\lesssim e^{-\left( a-\sqrt{\max\{a^{2}-m^{2}, 0\}}\right) t}\left( \|u_{0}\|_{H^{2k+j\sigma}}+\|u_{1}\|_{H^{\max\{2k+(j-1)\sigma,0\}}}\right) , \label{2.1} 
	\end{align}
\end{lemma}
Thanks to the exponential decay estimates, one can see here that the approach does not require any additional regularity for the initial data as happened for damped wave equations. So the energy space where the initial data lies is sufficient to get the desired estimates. The following integral inequality is proved in \cite{D'AbbiccoEbert}.
\begin{lemma}\cite[Lemma 6.3]{D'AbbiccoEbert}\label{second lemma}
	Let $c>0$ and $\alpha \in \mathbb{R}$. Then it holds
	$$\int_{0}^{t}e^{-c(t-s)}(1+s)^{-\alpha}ds\lesssim  (1+t)^{-\alpha}.$$
\end{lemma}
Now, the direct application of Lemma (4.1) in \cite{Cui} and Lemma \ref{second lemma} gives the following result.
\begin{lemma}\label{third lemma}
	Let $c>0$, $\beta>0$ and $\gamma\in (0,1)$. Then we have:
	$$\int_{0}^{t}e^{-c(t-s)}\int_{0}^{\tau}(\tau-s)^{-\gamma}(1+s)^{-\beta}dsd\tau\lesssim \left\{ 
	\begin{matrix}
	(1+t)^{-\gamma} \hfill& {if }\ \ \beta>1,
	&\cr 
	\\
	(1+t)^{1-\gamma-\beta}\hfill& {if }\ \ \beta<1.
	\end{matrix}\right.$$
\end{lemma}
We also need to use the following special case of fractional Gagliardo-Nirenberg inequality. We recall it as follows (see \cite{PhamKainaneReissig}).
\begin{lemma} \cite[Corollary 41] {PhamKainaneReissig}\label{FGN}
	Let $1<q<\infty$ and $\sigma>0$. Then, the following fractional Gagliardo-Nirenberg inequality holds for all $y\in H^{\sigma}(\mathbb{R}^{n})$
	$$\|y\|_{L^{q}(\mathbb{R}^{n})}\lesssim \|(-\Delta)^{\sigma/2}y\|_{L^{2}(\mathbb{R}^{n})}^{\theta_{q}}\,\|y\|_{L^{2}(\mathbb{R}^{n})}^{1-\theta_{q}},$$
	where $$\theta_{q}=\frac{n}{\sigma}\left(\frac{1}{2}-\frac{1}{q}\right)\in\left[0,1\right].$$
\end{lemma}
We are now in a position to stat our main result.
\section{global existence results} \label{GlobalExistence}   
\begin{theorem}\label{globalexistence1}
	Let us consider the Cauchy problem (\ref{1.1}) with $\sigma\geq 1$, $p>1$ and $\gamma\in(0,1)$ such that \begin{equation}\label{3.1}
	 \left\lbrace  
	 \begin{matrix}
	 1< p \leqslant \frac{n}{n-2\sigma} \hfill& {if }\ \ 2\sigma<n,
	 &\cr
	 \\
	1<p \hfill & {if } \ \
	n \leqslant 2\sigma. 
	 &\cr
	 \end{matrix}\right.  
	 \end{equation}
	 Moreover, we suppose the exponent $p$ satisfies 
	 \begin{equation}\label{3.2}  
	 p>\frac{1}{\gamma}.
	 \end{equation}
	Then, there exists a constant $\varepsilon>0$ such that for any data 
    $$u_{0} \in H^{\sigma}(\mathbb{R}^{n}), \ \ u_{1} \in L^{2}(\mathbb{R}^{n}) \ \ with\ \  \left\|(u_{0}, u_{1})\right\|_{H^{\sigma}(\mathbb{R}^{n})\times L^{2}(\mathbb{R}^{n})}<\varepsilon,$$ 
we have a unique global (in time) solution
$$u\in\mathcal{C}\left([0,\infty), H^{\sigma}(\mathbb{R}^{n})\right)\cap \mathcal{C}^{1}\left([0,\infty),L^{2}(\mathbb{R}^{n})\right)$$
to (\ref{1.1}). Furthermore, the solution satisfies the estimates:
	\begin{equation*}
	\|u(t,\cdot)\|_{L^{2}}\lesssim (1+t)^{-\gamma}\left\|(u_{0}, u_{1})\right\|_{H^{\sigma}(\mathbb{R}^{n})\times L^{2}(\mathbb{R}^{n})}, 
	\end{equation*}
	\begin{equation*}
	\|(-\Delta)^{\sigma/2}u(t,\cdot)\|_{L^{2}}\lesssim (1+t)^{-\gamma}\left\|(u_{0}, u_{1})\right\|_{H^{\sigma}(\mathbb{R}^{n})\times L^{2}(\mathbb{R}^{n})},
	\end{equation*}
	\begin{equation*}
	\|\partial_{t}u(t,\cdot)\|_{L^{2}}\lesssim (1+t)^{-\gamma}\left\|(u_{0}, u_{1})\right\|_{H^{\sigma}(\mathbb{R}^{n})\times L^{2}(\mathbb{R}^{n})}.
	\end{equation*}
\end{theorem} 
\begin{remark}
For higher space dimension $n>2\sigma$, we remark that the range $$\left( \frac{1}{\gamma}, \frac{n}{n-2\sigma}\right] $$ of exponent $p$  is not empty if and only if $$1-\frac{2\sigma}{n}<\gamma<1.$$ 
\end{remark}
Now we are ready to prove our theorem. The method is now standard and is based on the contraction mapping principle. 
\begin{proof}
The linear equation in (\ref{1.1}) has constant coefficients. So, by applying Duhamel's principle to the Cauchy problem (\ref{1.1}) we can write the solution $u$ as follows :
\begin{equation*}
u(t,x)=u^{lin}(t,x)+u^{nol}(t,x)
\end{equation*}
where 
$$u^{lin}(t,x)=K_{1}(t,x)\ast u_{0}(x)+K_{2}(t,x)\ast u_{1}(x)$$
and 
$$u^{nol}(t,x)=\Int_{0}^{t}K_{2}(t-\tau,x)\ast\left( \int_{0}^{\tau} (\tau-s)^{-\gamma}|u(s,x)|^{p}\right) dsd\tau.$$
Here the kernels $K_{i}$ are defined in \cite{d'abbicco2015} for $i=1,2$, and $\ast$ denotes the convolution product with respect to $x$. Let us now define for $T>0$ the following Banach space where the solution lives:
$$
X(T):=\mathcal{C}\left([0,T], H^{\sigma}(\mathbb{R}^{n})\right)\cap \mathcal{C}^{1}\left([0,T],L^{2}(\mathbb{R}^{n})\right),
$$
we equip this space with a suitable norm as follows:
\begin{align}
\|u\|_{B(T)}&=\sup_{0\leq t\leq T}\Big((1+t)^{\gamma}\left( \|u(t,\cdot)\|_{L^{2}}+\|(-\Delta)^{\sigma/2}u(t,\cdot)\|_{L^{2}}+\|\partial_{t}u(t,\cdot)\|_{L^{2}}\right) \Big). \label{3.3}
\end{align}
We remark that the weight function $(1+t)^{\gamma}$ is chosen in an appropriate way to prove our result and not chosen from linear estimates as usually happened for power nonlinearity.

 Let us define the operator $\mathcal{S}$ on the Banach space $X(T)$ by:
\begin{equation*}
\mathcal{S} : X(T)\longrightarrow X(T) : u\longmapsto \mathcal{S}u=u^{lin}+u^{nol}.   
\end{equation*}
The main  steps is to prove the following two inequalities for the operator $\mathcal{S}$:
	\begin{align}
\|\mathcal{S}u\|_{X(T)} &\lesssim \left\|(u_{0},u_{1})\right\|_{H^{\sigma}(\mathbb{R}^{n})\times L^{2}(\mathbb{R}^{n})}+\|u\|_{X(T)}^{p},  \ \ \forall u\in X(T),\label{3.4} \\ 
\|\mathcal{S}u-\mathcal{S}v\|_{X(T)} &\lesssim \|u-v\|_{X(T)}\Big( \|u\|_{X(T)}^{p-1}+\|v\|_{X(T)}^{p-1}\Big),  \ \ \forall (u,v) \in X(T). \label{3.5} 
\end{align}
\textit{Step 1:} We begin with the first inequality (\ref{3.4}). The proof is equivalent to 
$$\|u^{lin}\|_{X(T)}\lesssim \left\|(u_{0},u_{1})\right\|_{H^{\sigma}(\mathbb{R}^{n})\times L^{2}(\mathbb{R}^{n})}  \ \ \  and \ \ \ \|u^{nol}\|_{X(T)}\lesssim \|u\|^{p}_{X(T)}.$$ 
The first inequality trivially satisfied. Because we have:
\begin{align*}
\|u^{lin}\|_{X(T)}&=\sup_{0\leq t\leq T}\Big((1+t)^{\gamma}\Big( \|u^{lin}(t,\cdot)\|_{L^{2}}+\|(-\Delta)^{\sigma/2}u^{lin}(t,\cdot)\|_{L^{2}}+\|\partial_{t}u^{lin}(t,\cdot)\|_{L^{2}}\Big) \Big),
\end{align*}
\begin{align}
&\lesssim\sup_{0\leq t\leq T}(1+t)^{\gamma}e^{-\left( a-\sqrt{\max\{a^{2}-m^{2}, 0\}}\right) t}\left\|(u_{0},u_{1})\right\|_{H^{\sigma}(\mathbb{R}^{n})\times L^{2}(\mathbb{R}^{n})} \nonumber \\
&\hspace{3cm} 
\lesssim \left\|(u_{0},u_{1})\right\|_{H^{\sigma}(\mathbb{R}^{n})\times L^{2}(\mathbb{R}^{n})}.  \label{3.6}
\end{align}
Now we turn to second inequality. The proof needs to estimates the following norms:
$$\|u^{nol}(t,\cdot)\|_{L^{2}}, \ \ \ \|\partial_{t}u^{nol}(t,\cdot)\|_{L^{2}}, \ \ \ \|(-\Delta)^{\sigma/2}u^{nol}(t,\cdot)\|_{L^{2}}.$$
The above lemmas from Section \ref{preliminaries and tools} come into play to estimates these terms. We have: 
$$\|u^{nol}(t,\cdot)\|_{L^{2}}\lesssim \Int_{0}^{t}e^{-\left( a-\sqrt{\max\{a^{2}-m^{2}, 0\}}\right)(t-s)}\int_{0}^{\tau} (\tau-s)^{-\gamma}\|u(s,\cdot)\|_{L^{2p}}^{p} dsd\tau, $$
as well as
 $$\|\partial_{t}u^{nol}(t,\cdot)\|_{L^{2}}\lesssim \Int_{0}^{t}e^{-\left( a-\sqrt{\max\{a^{2}-m^{2}, 0\}}\right)(t-s)}\int_{0}^{\tau} (\tau-s)^{-\gamma}\|u(s,\cdot)\|_{L^{2p}}^{p} dsd\tau, $$
 
 $$\|(-\Delta)^{\sigma/2}u^{nol}(t,\cdot)\|_{L^{2}}\lesssim \Int_{0}^{t}e^{-\left( a-\sqrt{\max\{a^{2}-m^{2}, 0\}}\right)(t-s)}\int_{0}^{\tau} (\tau-s)^{-\gamma}\|u(s,\cdot)\|_{L^{2p}}^{p} dsd\tau. $$
In order to estimates the norm $\|u(s,\cdot)\|_{L^{2p}}^{p}$ we use the fractional Galgiardo-Nirenberg inequality from Lemma \ref{FGN} we have:
$$\|u(s,\cdot)\|_{L^{2p}}^{p}\lesssim \|(-\Delta)^{\sigma/2}u(s,\cdot)\|_{L^{2}}^{p\theta_{2p}}\|u(s,\cdot)\|_{L^{2}}^{p(1-\theta_{2p})}.$$
Here, we have the fact that:
$$\|(u(s,\cdot),\ (-\Delta)^{\sigma/2}u(s,\cdot))\|_{L^{2}}\lesssim (1+s)^{-\gamma}\|u\|_{X(T)},$$
this leads to 
$$\|u(s,\cdot)\|_{L^{2p}}^{p}\lesssim (1+s)^{-\gamma p}\|u\|_{X(T)}^{p}.$$
Putting this into all the integrals and we use $\gamma p>1$ from Lemma \ref{third lemma}, then we arrive to the following desired estimates:
$$\|u^{nol}(t,\cdot)\|_{L^{2}} \lesssim (1+t)^{-\gamma}\|u\|_{X(T)}^{p},$$
$$\|(-\Delta)^{\sigma/2}u^{nol}(t,\cdot)\|_{L^{2}}\lesssim (1+t)^{-\gamma}\|u\|_{X(T)}^{p},$$
$$\|\partial_{t}u^{nol}(t,\cdot)\|_{L^{2}}\lesssim (1+t)^{-\gamma}\|u\|_{X(T)}^{p},$$
in this way we complete the proof of (\ref{3.4}). The proof of (\ref{3.5}) can be done analogously. In fact, we choose two elements $(u, v) \in X(T)\times X(T)$ and we write:
$$\mathcal{S}u(t,x)-\mathcal{S}v(t,x)=\Int_{0}^{t}K_{2}(t-\tau,x)\ast\left( \int_{0}^{\tau} (\tau-s)^{-\gamma}(|u(s,x)|^{p}-|v(s,x)|^{p})\right) dsd\tau.$$
As above we have the following estimates :
$$\|\mathcal{S}u(t,\cdot)-\mathcal{S}v(t,\cdot)\|_{L^{2}}\lesssim \Int_{0}^{t}e^{-\left( a-\sqrt{\max\{a^{2}-m^{2}, 0\}}\right)(t-s)}\int_{0}^{\tau} (\tau-s)^{-\gamma}\||u(s,\cdot)|^{p}-|v(s,\cdot)|^{p}\|_{L^{2}} dsd\tau $$
as well as
$$\|\partial_{t}(\mathcal{S}u(t,\cdot)-\mathcal{S}v(t,\cdot))\|_{L^{2}}\lesssim \Int_{0}^{t}e^{-\left( a-\sqrt{\max\{a^{2}-m^{2}, 0\}}\right)(t-s)}\int_{0}^{\tau} (\tau-s)^{-\gamma}\||u(s,\cdot)|^{p}-|v(s,\cdot)|^{p}\|_{L^{2}}  dsd\tau, $$

$$\|(-\Delta)^{\sigma/2}(\mathcal{S}u(t,\cdot)-\mathcal{S}v(t,\cdot))\|_{L^{2}}\lesssim \Int_{0}^{t}e^{-\left( a-\sqrt{\max\{a^{2}-m^{2}, 0\}}\right)(t-s)}\int_{0}^{\tau} (\tau-s)^{-\gamma}\||u(s,\cdot)|^{p}-|v(s,\cdot)|^{p}\|_{L^{2}}  dsd\tau. $$
By employing the H\"{o}lder's inequality, we derive the following
\begin{equation*}
\| |u(s,\cdot)|^{p}-|v(s,\cdot)|^{p}\|_{L^{2}}\leq \|u(s,\cdot)-v(s,\cdot)\|_{L^{2 p}}\left(\|u(s,\cdot)\|_{L^{2 p}}^{p-1}+\|v(s,\cdot)\|_{L^{2p}}^{p-1} \right).
\end{equation*}
 Using again the definition of the norms $\|u\|_{X(T)}$, $\|u-v\|_{X(T)}$ and the fractional Gagliardo-Nirenberg inequality we have 
 \begin{equation*}
 \| |u(s,\cdot)|^{p}-|v(s,\cdot)|^{p}\|_{L^{2}}\lesssim (1+s)^{-\gamma p} \|u-v\|_{X(t)}\left(\|u\|_{X(T)}^{p-1}+\|v\|_{X(T)}^{p-1} \right).
 \end{equation*} 
 The same condition $p\gamma>1$ leads to the desired estimates. Summarizing, the proof of Theorem \ref{globalexistence1} is completed.
\end{proof}
\begin{remark}
It is also more interesting to study the interaction between the two equations of the form (\ref{1.1}), that is the following weakly coupled system :
\begin{equation}\label{3.7}
\left\lbrace 
\begin{array}{ll}
\partial_{t}^{2}u+(-\Delta)^{\sigma_{1}}u+2a_{1}\partial_{t}u+m_{1}^{2}u=\Int_{0}^{t}(t-s)^{-\gamma_{1}}\left|v(s,x)\right| ^{p}ds\\
\partial_{t}^{2}v+(-\Delta)^{\sigma_{2}}v+2a_{2}\partial_{t}v+m_{2}^{2}v=\Int_{0}^{t}(t-s)^{-\gamma_{2}}\left|u(s,x)\right| ^{q}ds\\
\\
u(0,x)=u_{0}(x),\ \ \partial_{t}u(0,x)=u_{1}(x), \ \ v(0,x)=v_{0}(x),\ \ \partial_{t}v(0,x)=v_{1}(x),
\end{array}\right.
\end{equation}
where 
\begin{equation*}
(t,x)\in \mathbb{R_{+}}\times \mathbb{R}^{n},, \ \ a_{1},a_{2}>0, \ \ m_{1},m_{2}>0, \ \ \gamma_{1}, \gamma_{2}\in(0,1), \ \ p, q >1 \ \ \ and \ \ \sigma_{1}, \sigma_{2}\geq 1.
\end{equation*} 
The goal is to find the $p-q$ curve which guarantees the global (in time) existence of small data solutions to (\ref{3.7}).
\end{remark}

\end{document}